\documentclass[12pt,thmsa]{article}
\usepackage{amsmath}
\usepackage{graphicx}
\usepackage{amsfonts}
\usepackage{amssymb}

\newtheorem{Cor}{Corollary}

\newtheorem{Lem}{Lemma}

\newtheorem{Assumpt}{Assumption}
\newtheorem{Theo}{Theorem}

\newcommand{\R}{\mathbb R}

\newcommand{\p}{\mathbb P}
\newcommand{\E}{\mathbb E}

\begin{document}

\begin{center}
{\Large Strong consistency of kernel  estimator in a semiparametric regression model}

{\Large \bigskip}

Emmanuel de dieu NKOU and Guy Martial \ NKIET

\bigskip

Universit\'{e} des Sciences et Techniques de Masuku

BP 943 \ Franceville, Gabon

E-mail : emmanueldedieunkou@gmail.com,  gnkiet@hotmail.com.

\bigskip
\end{center}

\noindent\textbf{Abstract.} Estimating the  effective dimension reduction (EDR) space, related to the semiparametric regression model introduced by Li \cite{sir}, is based on the estimation of the covariance matrix $\Lambda$ of  the conditional expectation of the vector of predictors given the response. An estimator $\widehat{\Lambda}_n$ of $\Lambda $ based on kernel method was introduced by Zhu and Fang \cite{Asymptotics} who then derived, under some conditions, the asymptotic distribution of $\sqrt{n}\left(\widehat{\Lambda}_n-\Lambda\right)$, as $n\rightarrow +\infty$. In this paper, we obtain, under specified conditions, the almost sure convergence of  $\widehat{\Lambda}_n$ to  $\Lambda$, as  $n\rightarrow +\infty$. 

\bigskip



\section{Introduction}

\label{intro}
Given a univariate response variable $Y$, we consider the regression model:
\begin{equation}\label{a1}
Y=F(\beta_1^{T}X,...,\beta_{N}^{T}X,\varepsilon), 
\end{equation}
where $X$ is a $d$-dimensional random vector with covariance matrix assumed, without loss of generality,  to be the identity matrix, $N$ is an integer of $\mathbb{N}^\ast$ such that $N<d$, $\beta_1,\cdots,\beta_N$ are vectors in $\mathbb{R}^d$, and $\varepsilon$ is a real random variable that is independent  of $X$, and $F$ is an arbitrary unknown function on $\mathbb{R}^{N+1}$. This model, introduced by Li  \cite{sir}, permits to achieve dimension reduction since the number $N$  of variables to be considered for estimating  $F$ is less than the initial dimension $d$ of the regressor vector $X$. It expresses the fact that the projection of $X$ onto the $N$-dimensional subspace spanned by   $\beta_1^{T}X,...,\beta_{N}^{T}X$, named the effective dimension reduction (EDR) space,  contains all information about the response variable $Y$. Estimating $N$ and the EDR space  is then a crucial issue that has  been tackled in several works (e.g.,    [1,2,3,4,8,9,10,11,12,13,16]). Since the directions $\beta_1,\cdots,\beta_K$ are, under some conditions,  characterized  as eigenvectors of the covariance matrix  $\Lambda $ of $\mathbb{E}(X\vert Y)$, the aforementioned estimation problem is based on estimation of $\Lambda$.  The most popular method for doing that is based on slicing the range of $Y$ and leads to the well known sliced inverse regression (SIR) method that was introduced by Li  \cite{sir} (see also \cite{Hsing}). An alternative method was  introduced by Zhu and Fang \cite{Asymptotics}; in this work an estimator $\widehat{\Lambda}_n$ of $\Lambda $  based on kernel method is proposed and the limiting distribution of $\sqrt{n}\left(\widehat{\Lambda}_n-\Lambda\right)$, as $n\rightarrow +\infty$,  is derived under some conditions.  Since this result just implies weak  consistency  of  $\widehat{\Lambda}_n$, that is the convergence in probability of $\widehat{\Lambda}_n$ to $\Lambda$ as $n\rightarrow +\infty$, it is natural to wonder if one could obtain strong consistency for  $\widehat{\Lambda}_n$.

In this paper, we tackle this problem and we prove, under some conditions, the almost sure convergence of  $\widehat{\Lambda}_n$ to $\Lambda$, as $n\rightarrow +\infty$. The paper is organized as follows: Section 2 is devoted to the presentation of the used estimator, that is the estimator given in \cite{Asymptotics}. In Section 3, the assumptions needed for our results are given, and the main theorems that establish the aforementioned consistency are given. Then, the  proofs of all lemmas and theorems are postponed in Section 4.
\section{Preliminaries and notations}
\label{notations}

Letting $f$ be the density of $Y$, we suppose that, for all $y\in\mathbb{R}$, we have  $f(y)>0$; then,    for any $j=1,\cdots,d$, we consider
\[
R_j(y)=\mathbb{E}(X_j\vert Y=y)=\frac{g_j(y)}{f(y)}
\,\,\textrm{ where }\,\,g_j(y)=\int_\mathbb{R} z f_{_{(X_j,Y)}}(z,y)dz,
\] 
$f_{_{(X_j,Y)}}$ being the density of the pair $(X_j,Y)$. Then, we consider the random vector
\[
R(Y)=\bigg(R_1(Y),...,R_d(Y)\bigg)^T=\bigg(\mathbb{E}(X_1\vert Y),...,\mathbb{E}(X_d\vert Y)\bigg)^T=\mathbb{E}\bigg(X\vert Y\bigg)
\]
and its covariance matrix $
\Lambda = Cov\bigg(\mathbb{E}\left(X\vert Y\right)\bigg) $ which  is of great importance since the EDR space is obtained from its spectral analysis (e.g.  \cite{sir}). It cannot be computed in practice since it  depends  on the distribution of $(X,Y)$ which is generally unknown; that is why
approaches for its estimation have been investigated by several authors. Li \cite{sir}  considered an estimation method based on slicing the range of $Y$, so introducing sliced inverse regression, whereas Zhu and Fang \cite{Asymptotics} introduced a kernel estimator. More precisely, considering an i.i.d. sample $(Y_i,\,X_i)_{i=1,...n}$ of the pair   $\left(Y,\,X\right)$ of random variables
connected according  to model (\ref{a1}) and putting
\[
X_i=\left(X_{i1},\cdots,X_{id}\right)^T,
\]
we define kernel estimates of $f$ and the $g_j$'s  by:
\begin{equation*}\label{e5}
\widehat{f}_n(y)=\frac{1}{n}\sum_{i=1}^n \frac{1}{h_n}K\left(\frac{y-Y_i}{h_n}\right),\quad  \widehat{g}_{j,n}(y)=\frac{1}{n}\sum_{i=1}^n X_{ij}\frac{1}{h_n}K\left(\frac{y-Y_i}{h_n}\right), 
\end{equation*} 
where $h_n$ is a bandwidth and $K(\cdot)$ is a kernel function.  In order to avoid  small values in the denominator, Zhu and Fang (1996) proposed   to consider
\[
f_{b_n}(y)=\max\big(f(y),b_n\big)\,\,\textrm{ and }\,\, \widehat{f}_{b_n}(y)=\max\left(\widehat{f}_n(y),b_n\right),
\]
where  $\left(b_n\right)_{n\in\mathbb{N}^\ast}$ is a sequence of positive real numbers that satisfies the property:  $\lim_{n\rightarrow +\infty}(b_n)=0$.
Then, the $R_{b_n,j}$'s  defined  by 
\[
R_{b_n,j}(y)=\frac{g_j(y)}{f_{b_n}(y)}
\]
are estimated by 
\begin{equation*}
\widehat{R}_{b_n,j}(y)=\frac{\widehat{g}_{j,n}(y)}{\widehat{f}_{b_n}(y)}
\end{equation*}
and putting
\[
\widehat{R}_{b_n}(y)=\left(\widehat{R}_{b_n,1}(y),...,\widehat{R}_{b_n,d}(y)\right)^T,
\]
we take as estimator of   $\Lambda$ the random matrix: 
 \begin{equation*}
\widehat{\Lambda}_n=\frac{1}{n}\sum_{i=1}^{n}\left(\widehat{R}_{b_n}(Y_i)\right)\left(\widehat{R}_{b_n}(Y_i)\right)^T.
\end{equation*}
This estimator was considered in  \cite{Asymptotics} who then proved that $\sqrt{n}\left(\widehat{\Lambda}_n-\Lambda\right)$ converges in distribution, as $n\rightarrow +\infty$, to a normal distribution. This result implies that 
$\widehat{\Lambda}_n$ converges in probability, as $n\rightarrow +\infty$, to $\Lambda$, that is weak consistency of the estimator.  
In the following section, we establish, under specified conditions,    the almost sure convergence of  $\widehat{\Lambda}_n$ to $\Lambda$ as $n\rightarrow +\infty$.

\section{Assumptions and main results}

In this section, we present our assumptions, then we give the main results that establish almost sure convergence of  $\widehat{\Lambda}_n$ to $\Lambda$ as $n\rightarrow +\infty$.
\label{assumptions}

\begin{Assumpt}\label{ar2}
The random variable $X$ is bounded, i.e.  there exists   $G>0$ such that
$\Vert X\Vert_d\leq G$, where $\Vert\cdot\Vert_d$ is the usual Euclidean norm of $\mathbb{R}^d$.
\end{Assumpt}

\begin{Assumpt}\label{ar12}
The random variable $Y$ has a bounded density $f$.
\end{Assumpt}

\begin{Assumpt}\label{ar3}The $g_j$'s and 
 $f$ are $3$-times differentiable and their  third derivatives  satisfy  the following condition: there exists a
neighborhood of the origin, say  $U$, and a constant $c>0$ such that, for any  $u\in U$, 
\[
\left|f^{(3)}\left(y+u\right)-f^{(3)}\left(u\right)\right|\leq c|u|\,\,\,\textrm{  and  }\,\,\,\left|g_j^{(3)}\left(y+u\right)-g_j^{(3)}\left(u\right)\right|\leq c|u|,
\]
for $j=1,\cdots,d$.
\end{Assumpt}

\begin{Assumpt}\label{reg1}
For any pair $(k,\ell)$ such that  $1\leq k,\ell\leq d$, and any   $u\in U$,
 \[
\left|R_k(y+u)R_\ell(y+u)-R_k(y)R_\ell(y)\right|\leq c|u|.
\]
\end{Assumpt}

\begin{Assumpt}\label{reg2}
There exists an integer $r>6$ such that, for any $j\in \left\{1,\cdots,d\right\}$, the function $ g_j$ belongs to the set
\[
 \Sigma\left(r,L,\alpha\right)=\left\{g\in \mathcal{D}^r \quad /\quad \forall \left(x,\,y\right), \,\left|g^{(r)}(x)- g^{(r)}(y) \right|\leq L\left|x-y\right|^{\alpha}\right\},
\]
 where $\alpha \in ]0,1]$, $L>0$  and  $\beta := r+\alpha$ satisfies  $\beta> 7$, and $\mathcal{D}^r$ denotes the space of $r$-times differentiable functions.
 
\end{Assumpt}

\begin{Assumpt}\label{ar1}

\begin{enumerate}
  \item[(i)] The kernel $K$ is continuous and its   support   is the interval $[-1,\,1]$;
  \item[(ii)] $K$ is symmetric about $0$;
  \item[(iii)] The kernel $K$ is bounded, that is: 
$
\sup_{u\in\mathbb{R}} \left|K(u)\right|=D<+\infty.       
$
  \item[(iv)] The kernel $K$ is of order $r$, that is  
\[\int  u^k K(u)du=0\,\,\,\textrm{for }k\in \left\{1,2,\cdots,r\right\};
\] 
  \item[(v)] 
\[
\int  \left|K(u)\right|du< +\infty \,\,\,\textrm{ and }\,\,\,\int\left|u\right|^{\beta} \left|K(u)\right|du<+ \infty.
\]
\end{enumerate}
\end{Assumpt}

\begin{Assumpt}\label{e17}
 When $n$ is large enough $h_n\sim n^{-c_1}$ and $b_n\sim n^{-c_2}$ where $c_1$ and $c_2$ are  numbers satisfying  $c_1>0$, $0<c_2<1/10$  and  $1/8+c_2/4<c_1<1/4-c_2$.
\end{Assumpt}

\begin{Assumpt}\label{ar8}
The eigenvalues $\lambda_1,\cdots,\lambda_ d$ of $\Lambda$ verify:   $\lambda_1>\cdots>\lambda_ d>0$.
\end{Assumpt}

\bigskip

\noindent The  assumptions 3, 4, 6-$(i)$,  6-$(ii)$ and 7  was introduced in \cite{Asymptotics} and are necessary here to use some results of this paper. Assumption 2 concerns the density of $Y$ and is classical since it is satisfied for  the usual probability distributions.
The assumptions 5, 6-$(iv)$ and 6-$(v)$ are classical assumptions of nonparametric statistics literature (see, e.g.,  \cite{tsyb}). Assumption 6-$(iii)$ is satisfied, for  instance, by the gaussian kernel.

\bigskip

\noindent\textbf{Remark 1.} For overcoming technical difficulties  due to small values in the denominator, Zhu and Fang (1996) introduced the modified version $\widehat{f}_{b_n}=\max(\widehat{f}_n,b_n)$ of the  kernel estimate $\widehat{f}_n$ of the density $f$. But this approach does not  guarantee that we get a good estimator of $f$. Indeed,  if  we take $b_n=n^{-1/11}$, then $b_n$ is still larger than $1/2$ for very large values of $n$ (for example $n=2000$). So,  every value of $\widehat{f}_n$ could be  cut off and, therefore,   $\widehat{f}_{b_n}$ would have a constant value.   This is an undesirable property that makes $\widehat{f}_{b_n}$ a bad estimator of the density.  To overcome this problem, we can take $b_n=\min(a,n^{-c_2})$, where $a$ is a fixed strictly positive  number. When $a$ is sufficiently small $\widehat{f}_{b_n}$ is near from $\widehat{f}_{n}$  and is, therefore, a good estimate of $f$. Indeed, it is easy to check that 
$
\sup_{x\in\mathbb{R}}\vert \widehat{f}_{b_n}(x)-\widehat{f}_{n}(x)\vert\leq a$.
 Finally, by taking   $b_n=\min(a,n^{-c_2})$, we obtain a good estimate of the density and we still have $b_n\sim n^{-c_2}$ as required in Assumption 7.

\bigskip

\noindent  For a symmetric $(d\times d)$ matrix $A=(a_{k,\ell})_{1\leq k,\ell\leq d}$ , we denote by  $Vech(A)$ the  $d(d+1)/2$-dimensional  vector 
\[
(a_{11},\cdots,a_{d1},a_{22},a_{32},\cdots,a_{d2},a_{33},a_{43},\cdots,a_{d3},\cdots,a_{dd})^T.
\]   
For a vector $V=\left(v_1,v_2,\cdots,v_m\right)\in\mathbb{R}^m$, we denote $\left\|V\right\|_{\infty}=\max_{1\leq i\leq m}\left|v_i\right|$.

\bigskip

\noindent Now, we give results which establish  strong consistency for $\widehat{\Lambda}_n$ as estimator of $\Lambda$.
\bigskip

\begin{Theo}\label{e14}
Under the  assumptions \ref{ar2}, \ref{ar1} and \ref{e17} we have
\begin{equation*}
 \left\|Vech\Big(\widehat{\Lambda}_n - \mathbb{E}\left(\widehat{\Lambda}_n\right)\Big)\right\|_{\infty}= \mathbf{O}_{a.s.}\left(\left(\frac{\log n}{n}\right)^{\nu}\right)\quad 
\end{equation*}
with $\nu = 1/2-2(c_1+c_2)$.

\end{Theo}
\bigskip
Putting 
\[
\Lambda = \left(\lambda_{k,\ell}\right)_{1\leq k,\ell\leq d}\textrm{  and }\widehat{\Lambda}_n = \left(\widehat{\lambda}^{(n)}_{k,\ell}\right)_{1\leq k,\ell\leq d},
\]
we have:
\bigskip
\begin{Theo}\label{e13}
Under the assumptions   \ref{ar2} to  \ref{e17}, we have
for any $1\leq k,\ell\leq d$: 
\begin{equation*}\label{e26}
\lim_{n\rightarrow +\infty}\mathbb{E}\left(\widehat{\lambda}^{(n)}_{k,\ell}\right) =\lambda_{k,\ell}.  
\end{equation*}
\end{Theo}
\bigskip

\noindent The following theorem is our main result; it results from Theorems \ref{e14} and \ref{e13}.  
\bigskip

\begin{Theo}\label{e29}
Under the assumptions   \ref{ar2} to  \ref{e17},
$\widehat{\Lambda}_n $ converges almost surely to $\Lambda$, as $n\rightarrow +\infty$. 
\end{Theo}

\bigskip

\noindent As a consequence of this theorem, we can deduce strong consistency for estimators of the $\beta_k$'s. Since   the covariance matrix of $X$  is assumed to be equal to the $d\times d$ identity matrix $\mathbb{I}_d$ , then  $\beta_k$ (for $k=1,\cdots,N$) is an  eigenvector  of $\Lambda$  associated with the $k$-th largest eigenvalue $\lambda_k$  (see \cite{sir}). We consider the empirical covariance matrix 
\[
\widehat{\Sigma}_n=\frac{1}{n}\sum_{i=1}^n\left(X_i-\overline{X}_n\right)\left(X_i-\overline{X}_n\right)^T\,\,\,
\textrm{ where }\,\,\,
\overline{X}_n=\frac{1}{n}\sum_{i=1}X_i,
\]
and we denote by $\widehat{\eta}_k$ an eigenvector of $\widehat{\Lambda}_n$ associated with the $k$-th largest eigenvalue $\widehat{\lambda}_k$. Clearly, from strong law of large numbers,    $\widehat{\Sigma}_n$ converges almost surely to $\mathbb{I}_d$ as $n\rightarrow +\infty$; then  $\widehat{\Sigma}_n$ is also invertible for large values of $n$, and we can take as estimator of $\beta_k$ the vector $\widehat{\beta}_k=\widehat{\Sigma}_n^{-1/2}\widehat{\eta}_k$. Then, we have:

\bigskip

\begin{Cor}\label{cor1}
Under the assumptions   \ref{ar2} to  \ref{ar8}, for any $k\in\{1,\cdots,N\}$, 
$\widehat{\beta}_k $ converges almost surely to $\beta_k$, as $n\rightarrow +\infty$. 
\end{Cor}

\section{Proofs}
\label{proofs}
\subsection{Preliminary results}
First, we recall below a lemma given in \cite{Asymptotics}  (see p. 1058) and which will be useful for proving other results.
\begin{Lem}\label{ar15} Under the assumptions 1, 3, 4 and 7, we have almost surely: 
\[
\sup_{y\in\R}\left|\widehat{f}_n(y)-f(y)\right|=\mathbf{O}\left(h_n^4+n^{-1/2}h_n^{-1}\log n\right), 
\]
as $n\rightarrow +\infty$.
\end{Lem}

\bigskip

\begin{Lem}\label{l58} Under assumptions 1, 2,  3, 4 and 7, we have 
for any $j\in\{1,\cdots,d\}$: $\mathbb{E}\left(g_j^2(Y)\right)<+\infty$.
\end{Lem}
\textbf{Proof.}
According to  \cite{Asymptotics}  (see p. 1059), we have: 
\[
\sup_{y\in\R}\left|\E\left(\widehat{g}_{j,n}(y)\right) - g_j(y)\right|=\mathbf{O}(h_n^4).
\]
Then,  there exists $M_1>0$ such that  $
\sup_{y\in\R}\left|\E\left(\widehat{g}_{j,n}(y)\right) - g_j(y)\right|\leq M_1$
 for any $n\in\mathbb{N}^\ast$. On the other hand,  
\begin{eqnarray*}
\left\vert \mathbb{E}\left(\widehat{g}_{j,n}(y)\right)\right\vert&\leq &Gh_n^{-1}\mathbb{E}\left(\left\vert K\left(\frac{y-Y_1}{h_n}\right)\right\vert\right)=
Gh_n^{-1}\int  \left\vert K\left(\frac{y-t}{h_n}\right)\right\vert f(t)\, dt\\
&=&G\int \left\vert K(u)\right\vert f(h_ny-u)\,\,du\leq G\Vert f\Vert_\infty\int  \left\vert K(u)\right\vert du,
\end{eqnarray*}
where $\Vert f\Vert_\infty=\sup_{t\in\mathbb{R}}f(t)$. Therefore,  for any $y\in\mathbb{R}$, 
\begin{eqnarray*}
\vert g_j(y)\vert&\leq& \left|\E\left(\widehat{g}_{j,n}(y)\right) - g_j(y)\right|+\left|\E\left(\widehat{g}_{j,n}(y)\right)\right| \\
&\leq & M_1+ G\Vert f\Vert_\infty\int  \left\vert K(u)\right\vert du.
\end{eqnarray*}
This shows that $g_j(Y)$ is a bounded real random variable and, therefore,  $\mathbb{E}\left(g_j^2(Y)\right)<+\infty$.
\hfill $\Box$
\bigskip

\begin{Lem}\label{ar6}Under the  assumptions 5 and 6, we have
for any $y\in\mathbb{R}$ and any $j\in\{1,\cdots,d\}$:
\[
\left|\int g_j(y-uh_n)K(u)du -g_j(y)\right| \leq C\,h_n^{\beta},
\] 
where $C>0$.
\end{Lem}
\textbf{Proof.} By a Taylor expansion, we have: 
\[
g_j(y-uh_n)=g_j(y)\,+\,\sum_{k=1}^{r-1}\frac{g^{(k)}(y)}{k!}(-1)^k u^k h_n^k \,+\,\frac{(-1)^ru^rh_n^r}{r!} g_j^{(r)}(y-\theta u h_n),
\] 
where $\theta \in ]0,\,1[$. Thus,
\begin{eqnarray*}
\int g_j(y-uh_n)K(u)du &=&\int g_j(y)K(u)du\,+\,\sum_{k=1}^{r-1}\frac{g^{(k)}(y)}{k!}(-1)^k u^k h_n^k\int u^kK(u)du\\                        
&+&\frac{(-1)^rh_n^r}{r!} \int g_j^{(r)}(y-\theta u h_n)u^r K(u)du,\\
                       &=& g_j(y) \,+\,\frac{(-1)^rh_n^r}{r!} \int g_j^{(r)}(y-\theta u h_n)u^r K(u)du.
\end{eqnarray*}
Furthermore, since
\[
 \int g_j^{(r)}(y) u^rK(u)du=g_j^{(r)}(y)\int u^r K(u) du = 0,
\]
it follows:
\begin{eqnarray*}
\int g_j(y-uh_n)K(u)du - g_j(y) &=& \frac{(-1)^r h_n^r}{r!}\left\{\int g_j^{(r)}(y-\theta u h_n) u^r K(u)du - \int g_j^{(r)}(y)u^rK(u)du\right\}\\
                                &=& \frac{(-1)^r h_n^r}{r!}\int \left(g_j^{(r)}(y-\theta u h_n) - g_j^{(r)}(y)\right)u^r K(u)du.
\end{eqnarray*}
Thus, under Assumption 5,
\begin{eqnarray*}
\left|\int g_j(y-uh_n)K(u)du - g_j(y)\right| & \leq & \frac{ h_n^r}{r!}\int \left|g_j^{(r)}(y-\theta u h_n) - g_j^{(r)}(y)\right| \left|u\right|^r \left|K(u)\right|du\\
                                             &\leq & \frac{h_n^r}{r!}\int L\, \theta^{\alpha}\left|u\right|^{\alpha}  h_n^{\alpha}\left|u\right|^r \left|K(u)\right|du \\
                                             &\leq & h_n^{r+\alpha }\frac{L}{r!}\int \left|u\right|^{r+\alpha}\left|K(u)\right|du\\  
                                             &\leq &  h_n^{\beta}\frac{L}{r!}\int \left|u\right|^{\beta}\left|K(u)\right|du,
\end{eqnarray*}
what gives the required inequality with  $C=\frac{L}{r!}\int \left|u\right|^{\beta}\left|K(u)\right|du$. \hfill $\Box$

\begin{Lem}\label{ar7}
Considering
\[
\mathcal{E}_{j,n}=\int f(y) \left(\int g_j(y-uh_n) K(u) du\right)^2 dy\,\,\,\textrm{ and }\,\,\,\mathcal{E}_j =\int f(y)g_j^2(y) dy = \E\left(g_j^2(Y)\right),
\]
then, under the assumptions 1 to  7, we have :
$
\,\left|\mathcal{E}_{j,n} - \mathcal{E}_j \right|\leq C^2h_n^{2\beta} + 2C\,h_n^{\beta}\,\E\bigg(\left|g_j(Y)\right|\bigg).
$
\end{Lem}
\textbf{Proof.} Using the equality $a^2-b^2=\left(a-b\right)^2 + 2b\left(a-b\right)$, we obtain :
\begin{eqnarray*}
\mathcal{E}_{j,n} - \mathcal{E}_j &=& \int f(y)\left[\left(\int g_j(y-uh_n)K(u)du\right)^2 - \bigg(g_j(y)\bigg)^2\right]dy\\
                                  &=& \int f(y)\left[\left(\int g_j(y-uh_n)K(u)du - g_j(y)\right)^2 \right.\\
& &+ \left. 2\, g_j(y)\left(\int g_j(y-uh_n)K(u)du - g_j(y)\right)\right]dy.
\end{eqnarray*}
Thus
\begin{eqnarray*}
\left|\mathcal{E}_{j,n} - \mathcal{E}_j\right| &\leq & \int f(y)\left[\left(\int g_j(y-uh_n)K(u)du - g_j(y)\right)^2 \right.\\
& &+\left. 2\, \left|g_j(y)\right|\left|\left(\int g_j(y-uh_n)K(u)du - g_j(y)\right)\right|\right]dy.
\end{eqnarray*}
Then, from Lemma \ref{ar6}, it follows
\begin{eqnarray*}
\left|\mathcal{E}_{j,n} - \mathcal{E}_j\right| & \leq & \int f(y)\left[C^2 h_n^{2\beta} + 2 \left|g_j(y)\right|C\,h_n^{\beta}\right]dy\\
                                               &=& C^2 h_n^{2\beta}\int f(y)dy + 2\,Ch_n^{\beta}\int \left|g_j(y)\right|f(y)dy\\
 & = & C^2h_n^{2\beta} + 2\,Ch_n^{\beta}\E\left(\left|g_j(Y)\right|\right).
\end{eqnarray*}
\hfill $\Box$

\begin{Lem}\label{ar9}
Putting $\delta_n=nh_n\left[\left(1-\frac{1}{n}\right)\mathcal{E}_{j,n} -\mathcal{E}_j  \right]$, we have under the assumptions 1 to  7,  $\lim_{n\rightarrow +\infty} \delta_n = 0$.
\end{Lem}
\textbf{Proof.} First, 
\[
\left|\mathcal{E}_{j,n}\right| \leq \left|\mathcal{E}_{j,n}-\mathcal{E}_j\right| + \left|\mathcal{E}_j\right|\leq
 C^2h_n^{2\beta} + 2C\,h_n^{\beta}\,\E\bigg(\left|g_j(Y)\right|\bigg)+ \E\left(g_j^2(Y)\right).
\]
Therefore,
\begin{eqnarray*}
\vert\delta_n\vert&=&nh_n\left\vert\left(\mathcal{E}_{j,n} -\mathcal{E}_j\right) - \frac{1}{n}\mathcal{E}_{j,n} \right\vert \leq nh_n\left[\left|\mathcal{E}_{j,n} -\mathcal{E}_j\right| + \frac{1}{n}\left|\mathcal{E}_{j,n}\right| \right]\\
&\leq & nh_n^{\beta+1}\left[C^2h_n^{\beta} + 2C\E\bigg(\left|g_j(Y)\right|\bigg)\right] \\
& &+ \frac{1}{n}\left[\E\left(g_j^2(Y)\right) +C^2h_n^{2\beta} + 2C\,h_n^{\beta}\,\E\bigg(\left|g_j(Y)\right|\bigg)\right].
\end{eqnarray*}
Clearly,
\[
\lim_{n\rightarrow  +\infty}\bigg(  \frac{1}{n}\left[\E\left(g_j^2(Y)\right) +C^2h_n^{2\beta} + 2C\,h_n^{\beta}\,\E\bigg(\left|g_j(Y)\right|\bigg)\right] \bigg)= 0.
\]
On the other hand, since $h_n\sim n^{-c_1}$, it follows that $nh_n^{\beta+1}\sim n^{1-(\beta+1)c_1}$. Further, from  $\beta > 7$ we deduce that $\frac{1}{\beta+1} < \frac{1}{8} < c_1$, that is   $1-(\beta+1)c_1<0$. Thus $\lim_{n\rightarrow +\infty} nh_n^{\beta+1} =0$ and, therefore, 
\[
\lim_{n\rightarrow  +\infty} \bigg( nh_n^{\beta+1}\left[C^2h_n^{\beta} + 2C\,\E\bigg(\left|g_j(Y)\right|\bigg)\right]\bigg)=0.
\]
 Finally, $\lim_{n\rightarrow +\infty} \delta_n = 0$. 
\hfill $\Box$

\begin{Lem}\label{ar11}Under the assumptions 1 to  7, 
we have: 
\[
\E\left(\left(\widehat{g}_{j,n}(Y) - g_j(Y)\right)^2\right) =\mathbf{O}\left(\frac{1}{n\,h_n}\right).
\]
\end{Lem}
\textbf{Proof.}  Considering the  random variable $W_{i,j,n}=X_{ij}\,K\left(\frac{Y-Y_i}{h_n}\right)$, we have:
\begin{eqnarray}\label{rel}
\E\left(\widehat{g}_{j,n}(Y)^2\right)=\frac{1}{nh_n^2}\E\left(W_{1,j,n}^2\right) + \frac{1}{h_n^2}\left(1-\frac{1}{n}\right)\E\left(W_{1,j,n}W_{2,j,n}\right).
\end{eqnarray}
Clearly,  $\E\left(W_{1,j,n}^2\right)=h_nJ_n$, where $J_n=\int\int\frac{1}{h_n}K^2\left(\frac{u}{h_n}\right)V(y-u)f(u)dudy$ with
$V(y)=\int z^2\, f_{(X,Y)}(z,y)\,dz$, and from Theorem 2.1.1 in  \cite{Rao1} it is known that $\lim_{n\rightarrow \infty}J_n = J:=\int\int V(y)K^2(u)f(u)du\,dy$. Furthermore,
\begin{eqnarray*}
 \E\left(W_{1,j,n}\,W_{2,j,n}\right) &=& \int\left[\left(\int K\left(\frac{u}{h_n}\right)g_j(y-u)du\right)\left( \int K\left(\frac{v}{h_n}\right)g_j(y-v)dv  \right)\right]f(y)dy\\
                                     &=& \int f(y)\left[\int K\left(\frac{u}{h_n}\right)g_j(y-u)du\right]^2dy.
\end{eqnarray*}
Putting  $t=\frac{u}{h_n}$, we obtain
\[
\E\left(W_{1,j,n}\,W_{2,j,n}\right) = h_n^2\int f(y)\left[\int K(t)g_j(y-th_n)dt\right]^2dy = h_n^2\mathcal{E}_{j,n}.
\]
On the other hand $\E\left(g_j(Y)^2\right) = \mathcal{E}_j$. Then, we deduce from (\ref{rel}) that
\begin{eqnarray*}
\E\left(\widehat{g}_{j,n}(Y)^2\right) - \E\left(g_j(Y)^2\right) &=& \frac{1}{nh_n^2}\,h_n\,J_n + \frac{1}{h_n^2}\left(1-\frac{1}{n}\right)h_n^2\mathcal{E}_{j,n} - \mathcal{E}_j
                                                                = \frac{1}{nh_n}\left(J_n + \delta_n\right).
\end{eqnarray*}
 On the other hand, it is easy to check  that 
\begin{eqnarray*}
\E\left(\left[\widehat{g}_{j,n}(Y) - g_j(Y)\right]^2\right) &=& \E\left(\widehat{g}_{j,n}(Y)^2\right) - \E\left(g_j(Y)^2\right)\nonumber\\
& & -2 \left\{\E\bigg(\widehat{g}_{j,n}(Y)\,g_j(Y)\bigg) - \E\left(g_j(Y)^2\right)\right\}\nonumber\\
                                                           &=& \frac{1}{nh_n}\left(J_n + \delta_n\right) -2\Delta_{j,n},\label{ar10}
\end{eqnarray*}
where $\Delta_{j,n} = \E\bigg(\widehat{g}_{j,n}(Y)\,g_j(Y)\bigg) - \E\left(g_j(Y)^2\right)$. We have:
\begin{eqnarray*} 
\E\bigg(\widehat{g}_{j,n}(Y)g_j(Y)\bigg)&=& \E\left[g_j(Y)\frac{1}{nh_n}\sum_{i=1}^n X_{ij}K\left(\frac{Y-Y_i}{h_n}\right)\right]\\
                                         &=& \frac{1}{h_n}\E\left[g_j(Y)X_{1j}K\left(\frac{Y-Y_1}{h_n}\right)\right]\\
                                         &=& \frac{1}{h_n}\int_{(3)}z\,g_j(y)K\left(\frac{y-u}{h_n}\right)f_{_{(X_j,Y_{_1},Y)}}(z,u,y)\,dz\,du\,dy \\
                                         &=& \frac{1}{h_n}\int_{(3)}z\,g_j(y)K\left(\frac{y-u}{h_n}\right)f_{_{(X_j,Y_{_1})}}(z,u)f_{_Y}(y)\,dz\,du\,dy \\
                                         &=& \frac{1}{h_n} \int g_j(y)f(y)\left(\int K\left(\frac{y-u}{h_n}\right)\left(\int zf_{_{(X_j,Y_{_1})}}(z,u) \right)\,du\right)dy\\
                                         &=& \frac{1}{h_n} \int g_j(y)f(y)\left(\int K\left(\frac{y-u}{h_n}\right)g_j(u)\,du\right)dy.
\end{eqnarray*}
Putting  $t=\frac{y-u}{h_n}$, we  obtain  
\[
\E\bigg(g_j(Y)\widehat{g}_{j,n}(Y)\bigg) = \int g_j(y)f(y)\left(\int g_j(y-th_n)K(t)\,dt\right)dy.
\]
Hence
\begin{eqnarray*}
\Delta_{j,n} &=& \int g_j(y)f(y)\left(\int g_j(y-th_n)K(t)dt\right)dy - \int g_j^2(y)f(y);\\
             &=& \int g_j(y)f(y)\left[\int g_j(y-th_n)K(t)dt - g_j(y)\right]dy
\end{eqnarray*}
and
\begin{eqnarray*}
\left|\Delta_{j,n}\right| &\leq &\int \left|g_j(y)\right|f(y)\left[\left|\int g_j(y-th_n)K(t)dt - g_j(y)\right|\right]dy\\
&\leq& Ch_n^{\beta}\int \left|g_j(y)\right|f(y)dy = Ch_n^{\beta}\E\left(\left|g_j(Y)\right|\right),
\end{eqnarray*}
the second inequality coming from Lemma \ref{ar6}. Then, we have:
\begin{eqnarray*} 
E\left(\left[\widehat{g}_{j,n}(Y) - g_j(Y)\right]^2\right) &\leq& \frac{1}{nh_n}\left|J_n + \delta_n\right| + 2\left|\Delta_{j,n}\right|\\
                                                           &\leq& \frac{1}{nh_n}\left\{\left|J_n + \delta_n\right| + 2\,Cnh_n^{\beta+1}\E\left(\left|g_j(Y)\right|\right)\right\}.
\end{eqnarray*}
Moreover,  $nh_n^{\beta+1} \sim n^{1-(\beta +1)}c_1$ and, since $\beta + 1> 8$ and $\frac{1}{8}<c_1$ we have the inequality  $1-(\beta + 1)c_1 < 0$  which implies $\lim_{n\rightarrow \infty} nh_n^{\beta+1} = 0$.
Consequently, 
\[
\lim_{n\rightarrow +\infty}\bigg( \left|J_n + \delta_n\right| + 2\,Cnh_n^{\beta+1}\E\left(\left|g_j(Y)\right|\right)\bigg)= \left|J\right|,
\] 
from what we deduce that   $\E\left(\left(\widehat{g}_{j,n}(Y) - g_j(Y)\right)^2\right)= \mathbf{O}\left(\frac{1}{nh_n}\right)$.
\hfill $\Box$
 
\begin{Lem}\label{ar22} 
Under the assumptions 1, 2, 3, 4 and  7, we have:
\begin{equation*}\label{l79_2}
 \lim_{n\rightarrow +\infty}\E\left(\left[\widehat{g}_{j,n}(Y) - g_j(Y)\frac{\widehat{f}_{b_n}(Y)}{f_{b_n}(Y)}\right]^2\right)= 0.
\end{equation*}
\end{Lem} 
\textbf{Proof.}
We have:
\begin{eqnarray}\label{ineg}
& & \E\left(\left[\widehat{g}_{j,n}(Y) - g_j(Y)\frac{\widehat{f}_{b_n}(Y)}{f_{b_n}(Y)}\right]^2\right) \nonumber\\
&=&\E\left(\left[\bigg(\widehat{g}_{j,n}(Y) - g_j(Y)\bigg)+g_j(Y)\bigg(1-\frac{\widehat{f}_{b_n}(Y)}{f_{b_n}(Y)}\bigg)\right]^2\right) \nonumber\\
&\leq &2\E\left(\bigg(\widehat{g}_{j,n}(Y) - g_j(Y)\bigg)^2\right)+2\E\left(g_j^2(Y)\bigg(1-\frac{\widehat{f}_{b_n}(Y)}{f_{b_n}(Y)}\bigg)^2\right).
\end{eqnarray}
Equation 4.4 in  \cite{Asymptotics} and Lemma \ref{ar15} allow to obtain, almost surely,  the inequality 
\[
\sup_{y\in\mathbb{R}}\vert \widehat{f}_{b_n}(y)-f_{b_n}(y)\vert\leq \sup_{y\in\mathbb{R}}\vert \widehat{f}_{n}(y)-f(y)\vert\leq M_2\,\left(h_n^4+n^{-1/2}h_n^{-1}\log n\right),
\]
where $M_2$ is a  positive constant. Therefore, almost surely,
\[
\bigg(1-\frac{\widehat{f}_{b_n}(Y)}{f_{b_n}(Y)}\bigg)^2=\bigg(\frac{\widehat{f}_{b_n}(Y)-f_{b_n}(Y)}{f_{b_n}(Y)}\bigg)^2\leq \bigg( M_2\,b_n^{-1}\left(h_n^4+n^{-1/2}h_n^{-1}\log n\right)\bigg)^2,
\]
and, consequently, 
\begin{eqnarray}\label{ineg2}
\E\left(g_j^2(Y)\bigg(1-\frac{\widehat{f}_{b_n}(Y)}{f_{b_n}(Y)}\bigg)^2\right)\leq  \bigg(M_2\,b_n^{-1}\left(h_n^4+n^{-1/2}h_n^{-1}\log n\right)\bigg)^2\E\left(g_j^2(Y)\right).
\end{eqnarray}
Clearly,  $b_n^{-1}h_n^4\sim n^{c_2-4c_1}$  and  $b_n^{-1}n^{-1/2}h_n^{-1}\sim n^{c_1+c_2-1/2}$ as $n\rightarrow +\infty$. Since, under assumption 4, we have $c_2-4c_1<0$ and $c_1+c_2-1/2<0$, it follows that $\lim_{n\rightarrow +\infty}\left(b_n^{-1}\left(h_n^4+n^{-1/2}h_n^{-1}\log n\right)\right)=0$. Then, from (\ref{ineg2}), (\ref{ineg}) and Lemma \ref{ar11}, we deduce the required result.
\hfill $\Box$

\subsection{Proof of Theorem \ref{e14}}
Since the class of functions
$$
\mathcal{H}_n=\bigg\{h_{(k,l)}:y \longmapsto h_{k,l}(y)=\frac{1}{n}\frac{\widehat{g}_{k,n}(y)\widehat{g}_{l,n}(y)}{\widehat{f}_{b_n}^2(y)}, 1\leq k,l\leq d \bigg\},
$$
is finite, we deduce from Lemma 3 in  \cite{Gine1} that it is a  \textit{Vapnik-\v{C}ervonenkis} (VC) class of functions with respect to the envelope  
\[
\mathbf{h}= \max\left\{\left|h_{k,l}\right|: h_{k,l} \in \mathcal{H}_n, \,1\leq k,l\leq d\right\}.
\]
The related  covering number $\mathcal{N}\left(\mathcal{H}_n,\left\|\cdot\right\|_{L^2(P)},\varepsilon \left\|\mathbf{h}\right\|_{L^2(P)}\right)$   satisfies,  for all $\varepsilon \in ]0,\,1[$ and for all probability measures $P$ on $\left(S, \mathcal{S}\right)$, 
\[
\mathcal{N}\left(\mathcal{H}_n,\left\|\cdot\right\|_{L^2(P)},\varepsilon \left\|\mathbf{h}\right\|_{L^2(P)}\right)\leq \bigg(\frac{A}{\varepsilon}\bigg)^\nu,
\]
where $A$ and $\nu$ are postive constant named the VC characteristics of $ \mathcal{H}_n$.  Assumptions \ref{ar2} and \ref{ar1}  imply    $\mathbf{h} \leq \frac{D^2G^2}{n h_n^2b_n^2}$, then we obtain, for all $h\in \mathcal{H}_n$,  
 $$
 \E\left(h(Y)\right)\leq \frac{D^2G^2}{h_n^2nb_n^2}\quad \mbox{ and } \quad \E\left(h^2(Y)\right)\leq \frac{D^4G^4}{h_n^4n^2b_n^4}.
 $$ 
 Taking   $\mu_n= \frac{D^2G^2}{nb_n^2h_n^2}$ and $\sigma_n^2=\frac{D^4G^4}{h_n^4n^2b_n^4}$, we can apply  Talagrand's   inequality (see  \cite{Talagrand1} and  Proposition 2.2 in  \cite{Gine2}):
there exist positive  constants $K_1$ and $K_2$, depending only on  $A$ and $\nu$, such that for all $t\geqslant K_1\Bigg[\mu_n \log \frac{A\mu_n}{\sigma} + \sqrt{n}\sigma \sqrt{\log \frac{A\mu_n}{\sigma_n}}\Bigg]$,
\begin{eqnarray*}
& &\p\Bigg\{\sup_{h\in \mathcal{H}_n}\left|\sum_{i=1}^n \bigg\{h\left(Y_i\right)-\E\bigg( h\left(Y\right)\bigg)\bigg\}\right| > t\Bigg\}\\
&\leq & K_2 \exp\left\{-\frac{1}{K_2}\frac{t}{\mu_n}\log \left(1+\frac{t\mu_n}{K_2 \left(\sqrt{n}\sigma_n + \mu_n \sqrt{\log \frac{A\mu_n}{\sigma_n}}\right)^2}\right)\right\},
\end{eqnarray*}
that is,
\begin{eqnarray}\label{talagrand}
& &\p\left\{\sup_{1\leq k,l\leq d}\left|\frac{1}{n}\sum_{i=1}^n  \frac{\widehat{g}_{k,n}\left(Y_i\right)\widehat{g}_{l,n}\left(Y_i\right)}{\widehat{f}_{b_n}^2\left(Y_i\right)}-\E\left( \frac{\widehat{g}_{k,n}\left(Y\right)\widehat{g}_{l,n}\left(Y\right)}{\widehat{f}_{b_n}^2\left(Y\right)}\right)\right| > t\right\}\nonumber\\
&\leq & K_2 \exp\left\{-\frac{1}{K_2}\frac{tnb_n^2h_n^2}{D^2G^2}\log \left(1+\frac{th_n^2b_n^2}{K_2\,D^2G^2\left(1+\sqrt{\log A}\right)^2} \right)\right\}.
\end{eqnarray}
Since
$h_n\sim n^{-c_1}$ and $b_n\sim n^{-c_2}$ , we have
 $\lim_{n\rightarrow \infty} \left(\frac{h_n}{n^{-c_1}}\right)=B_1$ and $\lim_{n\rightarrow \infty} \left(\frac{b_n}{n^{-c_2}}\right)=B_2$, where $B_1>0$  and $B_2>0$.
 Thus for  $\varepsilon$ such as  $0<\varepsilon<\min \left(1,\,B_1,\,B_2\right)$ and $n$ is large enough, we have 
   $$B_1-\varepsilon <\frac{h_n}{n^{-c_1}}<B_1+\varepsilon\quad\mbox{ and }\quad B_2-\varepsilon <\frac{b_n}{n^{-c_2}}<B_2+\varepsilon,\,$$ 
that is
   $$n^{-c_1}\left(B_1-\varepsilon\right) <h_n<n^{-c_1}\left(B_1+\varepsilon\right)\quad\mbox{ and }\quad n^{-c_2}\left(B_2-\varepsilon\right) <b_n< n^{-c_2}\left(B_2+\varepsilon\right).$$
Then, putting  $\delta = \left(B_1-\varepsilon\right)^2\left(B_2-\varepsilon\right)^2$, we deduce from (\ref{talagrand}) that
\begin{eqnarray}\label{talagrand2}
& &\p\left\{\sup_{1\leq k,l\leq d}\left|\frac{1}{n}\sum_{i=1}^n  \frac{\widehat{g}_{k,n}\left(Y_i\right)\widehat{g}_{l,n}\left(Y_i\right)}{\widehat{f}_{b_n}^2\left(Y_i\right)}-\E\left( \frac{\widehat{g}_{k,n}\left(Y\right)\widehat{g}_{l,n}\left(Y\right)}{\widehat{f}_{b_n}^2\left(Y\right)}\right)\right| > t\right\}\\
&\leq & K_2 \exp\left\{-\frac{1}{K_2}\frac{\delta\, t}{D^2G^2} n^{1-2(c_1+c_2)}\times \log\left( 1+\frac{\delta\, t}{K_2\,D^2G^2 \left(1+\sqrt{\log A}\right)^2}n^{-2(c_1+c_2)}\right)\right\}\nonumber
\end{eqnarray}
Let us put   $t_n=\left(\frac{\log n}{n}\right)^{1/2-2(c_1+c_2)}$; 
Assumption \ref{e17} implies $0<c_1+c_2<1/4$ and, consequently,  that  $\alpha =1/2-2(c_1+c_2)$ is  strictly positive. Then, $\lim_{n\rightarrow +\infty}(\log n)^\alpha=+\infty$ and, therefore, putting  $U\,=\,K_1\,G\,D\sqrt{\log A}$ we have for  $n$ large enough 
$$
\left(\log n\right)^{\alpha}  \geqslant   2U =  U(1+1) \geqslant   U\left(1+\sqrt{\frac{\log A}{n}}\right) =  U\left(\frac{\sqrt{n}+\sqrt{\log A}}{\sqrt{n}}\right),
$$
that is
\begin{eqnarray*}
\left(\frac{\log n}{n}\right)^{\frac{1}{2}-2\left(c_1+c_2\right)} & \geqslant & K_1\frac{D^2G^2}{nb_n^2 h_n^2}\times\sqrt{\log A}\left(\sqrt{n}+\sqrt{\log A}\right)
\end{eqnarray*}
what means that
\[
t_n\geqslant K_1\Bigg[\mu_n \log \frac{A\mu_n}{\sigma _n} + \sqrt{n}\sigma \sqrt{\log \frac{A\mu_n}{\sigma_n}}\Bigg].
\]
Then, (\ref{talagrand2}) can be applied to $t_n$ and we obtain
\begin{eqnarray*}\label{talagrand3}
\p\left\{\sup_{1\leq k,l\leq d}\left|\frac{1}{n}\sum_{i=1}^n  \frac{\widehat{g}_{k,n}\left(Y_i\right)\widehat{g}_{l,n}\left(Y_i\right)}{\widehat{f}_{b_n}^2\left(Y_i\right)}-\E\left( \frac{\widehat{g}_{k,n}\left(Y\right)\widehat{g}_{l,n}\left(Y\right)}{\widehat{f}_{b_n}^2\left(Y\right)}\right)\right| > t_n\right\}\leq v_n
\end{eqnarray*}
where
\[
v_n=K_2 \exp\left\{-\frac{1}{K_2}\frac{\delta\sqrt{n}(\log n)^\alpha}{D^2G^2}  \log\left( 1+\frac{\delta\, (\log n)^\alpha}{K_2\,D^2G^2 \sqrt{n}\left(1+\sqrt{\log A}\right)^2}\right)\right\}.
\]
Clearly, $v_n\sim w_n$ as $n\rightarrow +\infty$, where
\[
w_n=K_2 \exp\left\{-\frac{\delta^2(\log n)^{2\alpha}}{K_2^2D^4G^4\left(1+\sqrt{\log A}\right)^2} \right\},
\]
and since $\sum_{n=1}^{+\infty}w_n<+\infty$, we deduce that $\sum_{n=1}^{+\infty}v_n<+\infty$. Then from the above inequality it follows that
\begin{eqnarray*}\label{talagrand3}
\sum_{n=1}^{+\infty}\p\left\{\sup_{1\leq k,l\leq d}\left|\frac{1}{n}\sum_{i=1}^n  \frac{\widehat{g}_{k,n}\left(Y_i\right)\widehat{g}_{l,n}\left(Y_i\right)}{\widehat{f}_{b_n}^2\left(Y_i\right)}-\E\left( \frac{\widehat{g}_{k,n}\left(Y\right)\widehat{g}_{l,n}\left(Y\right)}{\widehat{f}_{b_n}^2\left(Y\right)}\right)\right| > t_n\right\}<+\infty,
\end{eqnarray*}
and the required result is obtained from   Borel Cantelli's lemma.
\hfill $\Box$
\subsection{Proof of Theorem \ref{e13}}
Let us consider
\[
R_{b_n,j}(y)=\frac{g_j(y)}{f_{b_n}(y)}, \,\,\,
I_{kl}^{(1)}(y)=\frac{g_k(y)g_l(y)}{f_{b_n}^2(y)}=R_{b_n,k}(y)R_{b_n,l}(y),
\]
\[
I_{kl}^{(2)}(y)=\frac{g_k(y)\widehat{g}_{l,n}(y)+g_l(y)\widehat{g}_{k,n}(y)}{f_{b_n}^2(y)}=\frac{R_{b_n,k}(y)\widehat{g}_{l,n}(y)}{f_{b_n}(y)} + \frac{R_{b_n,l}(y)\widehat{g}_{k,n}(y)}{f_{b_n}(y)},
\]
and
\[
I_{kl}^{(3)}(y)=2R_{b_n,k}(y)R_{b_n,l}(y)\frac{\widehat{f}_{b_n}(y)}{f_{b_n}(y)}.
\]
Denoting by  $\widehat{\lambda}^{(n)}_{k,l}$ the element at the $k$-th row and the $l$-th column of the matrix $\widehat{\Lambda}_n$  , it is known from   \cite{Asymptotics}  (see pp.  1059-1060) that
\begin{eqnarray}\label{i6}
\widehat{\lambda}^{(n)}_{k,l} = \frac{1}{n}\sum_{i=1}^n\left(\frac{\widehat{g}_{k,n}\left(Y_i\right)\widehat{g}_{l,n}\left(Y_i\right)}{\widehat{f}_{b_n}^2\left(Y_i\right)}\right)&=&\frac{1}{n}\sum_{i=1}^n\left\{I_{kl}^{(1)}\left(Y_i\right)+I_{kl}^{(2)}\left(Y_i\right)-I_{kl}^{(3)}\left(Y_i\right)\right\}\nonumber \\
& & - A_n + B_n+C_n-D_n,
\end{eqnarray}
where
\[
A_n=\frac{1}{n}\sum_{i=1}^n\left\{g_k\left(Y_i\right)\bigg(\widehat{g}_{l,n}\left(Y_i\right)-g_l\left(Y_i\right)\bigg) + g_l\left(Y_i\right)\bigg( \widehat{g}_{k,n}\left(Y_i\right)- g_k\left(Y_i\right)\bigg)\right\}\left( \frac{\widehat{f}_{b_n}^2\left(Y_i\right) - f_{b_n}^2\left(Y_i\right)}{\widehat{f}_{b_n}^2\left(Y_i\right) f_{b_n}^2\left(Y_i\right)}\right),
\]
\[
B_n=\frac{1}{n}\sum_{i=1}^n\frac{\left(\widehat{g}_{k,n}\left(Y_i\right) - g_k\left(Y_i\right)\right)\left(\widehat{g}_{l,n}\left(Y_i\right) - g_l\left(Y_i\right)\right)}{\widehat{f}_{b_n}^2\left(Y_i\right)},
\]
\[
C_n=\frac{1}{n}\sum_{i=1}^n R_{(b_n,k)}\left(Y_i\right)R_{(b_n,l)}\left(Y_i\right)\frac{\left(\widehat{f}_{b_n}^2\left(Y_i\right) - f_{b_n}^2\left(Y_i\right)\right)^2}{\widehat{f}_{b_n}^2\left(Y_i\right)f_{b_n}^2\left(Y_i\right)},
\]
and
\[
D_n=\frac{1}{n}\sum_{i=1}^n\bigg(\widehat{f}_{b_n}^2\left(Y_i\right)-f_{b_n}^2\left(Y_i\right)\bigg)^2\frac{R_{b_n,k}\left(Y_i\right)R_{b_n,l}\left(Y_i\right)}{f_{b_n}^2\left(Y_i\right)}.
\]
First, we will obtain the rates of convergence of  the sequences  $\E\left(A_n\right)$, $\E\left(B_n\right)$, $\E\left(C_n\right)$ and $\E\left(D_n\right)$   to $0$ as $n\rightarrow +\infty$. Clearly, $
\E\left(A_n\right)=\E\left(A_n^{(k,l)}\right) + \E\left(A_n^{(l,k)}\right)$, where 
\[
A_n^{(k,l)} = g_k\left(Y\right)\bigg(\widehat{g}_{l,n}\left(Y\right)-g_l\left(Y\right)\bigg)\left( \frac{\widehat{f}_{b_n}^2\left(Y\right) - f_{b_n}^2\left(Y\right)}{\widehat{f}_{b_n}^2\left(Y\right) f_{b_n}^2\left(Y\right)}\right).
\] 
Further,
\begin{eqnarray}\label{ean}
\left|\E\left(A_n^{(k,l)}\right)\right| & \leq &  \frac{1}{b^4_n}\,\E\left\{\left|g_k\left(Y\right)\bigg(\widehat{g}_{l,n}\left(Y\right)-g_l\left(Y\right)\bigg) \left(\widehat{f}_{b_n}^2\left(Y\right) - f_{b_n}^2\left(Y\right)\right)\right|\right\}\nonumber\\
                         & \leq &  \frac{1}{b_n^4}\,\E\left\{\left|g_k\left(Y\right)\bigg(\widehat{g}_{l,n}\left(Y\right)-g_l\left(Y\right)\bigg)\right| \left|\bigg(\widehat{f}_{b_n}\left(Y\right) - f_{b_n}\left(Y\right)\bigg)^2\,\right.\right.\nonumber\\
& &+\left.\left.\,2f_{b_n}\left(Y\right)\bigg(\widehat{f}_{b_n}\left(Y\right) - f_{b_n}\left(Y\right)\bigg)\right|\right\}\nonumber\\
                         & \leq &  \frac{1}{b_n^4}\,\E\left\{\left|g_k\left(Y\right)\bigg(\widehat{g}_{l,n}\left(Y\right)-g_l\left(Y\right)\bigg)\right|\,\, \left(\widehat{f}_{b_n}\left(Y\right) - f_{b_n}\left(Y\right)\right)^2\right\}\nonumber\\
                  & & +\quad\frac{2}{b_n^4}\,\E\left\{\left|g_k\left(Y\right)\bigg(\widehat{g}_{l,n}\left(Y\right)-g_l\left(Y\right)\bigg)\right| \,\, f_{b_n}\left(Y\right)\left|\widehat{f}_{b_n}\left(Y\right) - f_{b_n}\left(Y\right)\right|\right\}\nonumber\\
&\leq &  \frac{1}{b_n^4}\,\E\left\{\left|g_k\left(Y\right)\bigg(\widehat{g}_{l,n}\left(Y\right)-g_l\left(Y\right)\bigg)\right|\,\, \left(\widehat{f}_n\left(Y\right) - f\left(Y\right)\right)^2\right\}\\
&  &   +\quad\frac{2}{b_n^4} \E\left\{\left|g_k\left(Y\right)\bigg(\widehat{g}_{l,n}\left(Y\right)-g_l\left(Y\right)\bigg)\right|\,\, f_{b_n}\left(Y\right)\left|\widehat{f}_n\left(Y\right) - f\left(Y\right)\right|\right\}\nonumber.
\end{eqnarray}
Putting $\alpha_n=h_n^4+n^{-1/2}h_n^{-1}\log n$ and using Lemma \ref{ar15},  Cauchy-Schwartz inequality and Lemma \ref{ar11}, we obtain
\begin{eqnarray*}
& &\E\left\{\left|g_k\left(Y\right)\bigg(\widehat{g}_{l,n}\left(Y\right)-g_l\left(Y\right)\bigg) \left(\widehat{f}_{b_n}^2\left(Y\right) - f_{b_n}^2\left(Y\right)\right)\right|\right\}\\
&\leq & \alpha_n^2\sqrt{\E\left(\bigg(\widehat{g}_{l,n}\left(Y\right)-g_l\left(Y\right)\bigg)^2\right)}\sqrt{\E\left(g_k^2\left(Y\right)\right)}
\leq  M_3\,\,\alpha_n^2\lambda_n^{1/2}\sqrt{\E\left(g_k^2\left(Y\right)\right)},
\end{eqnarray*}
where $\lambda_n=n^{-1}h_n^{-1}$ and $M_3\,$ is a positive constant. On the other hand, since for $n$ large enough  $f_{b_n}(Y)\leq\Vert f\Vert_\infty$, it follows
\begin{eqnarray*}
& &\E\left\{\left|g_k\left(Y\right)\bigg(\widehat{g}_{l,n}\left(Y\right)-g_l\left(Y\right)\bigg)\right|\,\, f_{b_n}\left(Y\right)\left|\widehat{f}_n\left(Y\right) - f\left(Y\right)\right|\right\}\\
&\leq &\alpha_n\Vert f\Vert_\infty\sqrt{\E\left(\bigg(\widehat{g}_{l,n}\left(Y\right)-g_l\left(Y\right)\bigg)^2\right)}\sqrt{\E\left(g_k^2\left(Y\right)\right)}
\leq M_3\, \alpha_n\lambda_n^{1/2} \Vert f\Vert_\infty \sqrt{\E\left(g_k^2\left(Y\right)\right)}.
\end{eqnarray*}
Therefore, from (\ref{ean}) we deduce that 
\begin{equation}\label{ar17}
 \left|\E\left(A_n\right)\right| = \mathbf{O}\left(b_n^{-4}\alpha_n\lambda_n^{1/2}\right)=\mathbf{O}\left(\beta_n\right)
\end{equation}
where $\beta_n=b_n^{-4}\alpha_n\lambda_n^{1/2}$. In addition, 
\begin{eqnarray*}
\E\left(B_n\right) &=& \E\left(\widehat{f}_{b_n}^{-2}\left(Y\right)\left(\widehat{g}_{k,n}\left(Y\right) - g_k\left(Y\right)\right)\left(\widehat{g}_{l,n}\left(Y\right) - g_l\left(Y\right)\right)\right)\\
&\leq & b_n^{-2}\sqrt{ \E\bigg(\left(\widehat{g}_{k,n}\left(Y\right) - g_k\left(Y\right)\right)^2\bigg)}\sqrt{ \E\bigg(\left(\widehat{g}_{l,n}\left(Y\right) - g_l\left(Y\right)\right)^2\bigg)}\\
&\leq &M_3^2\,b_n^{-2}\lambda_n;
\end{eqnarray*}
thus
\begin{equation}\label{ar18}
\left|\E\left(B_n\right)\right| = \mathbf{O}\left({b_n}^{-2}\lambda_n\right).
\end{equation}
Next, we have
\begin{eqnarray*}
\left|\E\left(C_n\right)\right| &\leq &  b_n^{-4} \E\left(\left|R_{k}\left(Y\right)R_{l}\left(Y\right)\right|\left(\widehat{f}_{b_n}^2\left(Y\right) - f_{b_n}^2\left(Y\right)\right)^2\right)\\
& \leq & b_n^{-4}\E\left(\left|R_{k}\left(Y\right)R_{l}\left(Y\right)\right| \bigg(\widehat{f}_{b_n}\left(Y\right) - f_{b_n}\left(Y\right)\bigg)^4 \right)\\
                                & &  +\, 4\,b_n^{-4}\E\left(\left|R_{k}\left(Y\right)R_{l}\left(Y\right)\right|f_{b_n}\left(Y\right)\bigg(\widehat{f}_{b_n}\left(Y\right) - f_{b_n}\left(Y\right)\bigg)^3\right)\\
& &  +\,  4\,b_n^{-4}\E\left(\left|R_{k}\left(Y\right)R_{l}\left(Y\right)\right|f_{b_n}^2\left(Y\right)\bigg(\widehat{f}_{b_n}\left(Y\right) - f_{b_n}\left(Y\right)\bigg)^2\right)\\
&\leq& b_n^{-4}\left(\alpha_n^4+4\Vert f\Vert_\infty \alpha_n^3+4\Vert f\Vert_\infty^2 \alpha_n^2\right)\E\left(\left|R_{k}\left(Y\right)R_{l}\left(Y\right)\right|\right).
\end{eqnarray*}
Thus
\begin{equation}\label{ar19}
\left|\E\left(C_n\right)\right| =\mathbf{O}\left(\ b_n^{-4}\alpha_n^2 \right).
\end{equation}
Similarly, we have
\begin{eqnarray*}
\left|\E\left(D_n\right)\right| &\leq & b_n^{-2}\E\left[\left|\widehat{f}_{b_n}^2\left(Y\right)-f_{b_n}^2\left(Y\right)\right|^2\left|R_{(b_n,k)}\left(Y\right)R_{(b_n,l)}\left(Y\right)\right|\right]\\
 & =&                           b_n^{-2}\E\left[\left|R_k\left(Y\right)R_l\left(Y\right)\right|\bigg(\widehat{f}_{b_n}\left(Y\right) - f_{b_n}\left(Y\right)\bigg)^4\right]\\
                              & &+\,4\,b_n^{-2}\E\left[\left|R_k\left(Y\right)R_l\left(Y\right)\right| f_{b_n}\left(Y\right)\bigg(\widehat{f}_{b_n}\left(Y\right) - f_{b_n}\left(Y\right)\bigg)^3\right]\\
                              & &+\,4\,b_n^{-2}\E\left[\left|R_k\left(Y\right)R_l\left(Y\right)\right| f_{b_n}^2\left(Y\right)\bigg(\widehat{f}_{b_n}\left(Y\right) - f_{b_n}\left(Y\right)\bigg)^2 \right]\\
                              &\leq & b_n^{-2}\left(\alpha_n^4+4\Vert f\Vert_\infty \alpha_n^3+4\Vert f\Vert_\infty^2 \alpha_n^2\right)\E\left(\left|R_{k}\left(Y\right)R_{l}\left(Y\right)\right|\right),
\end{eqnarray*}
what implies
\begin{equation}\label{ar20}
\left|\E\left(D_n\right)\right| = \mathbf{O}\left(\alpha_n^2b_n^{-2}\right).
\end{equation} 
From     Eq. (\ref{ar17}) to Eq.  (\ref{ar20}), we obtain
\begin{equation*}\label{ar21}
\left|\E\left(-A_n+B_n+C_n-D_n\right)\right|= \mathbf{O}\left(n^{-(1-4c_1-2c_2)}\log n\right).
\end{equation*}
Then, from Eq.(\ref{i6}) we deduce that
\begin{equation}\label{i9}
\E\left(\widehat{\lambda}^{(n)}_{k,l}\right)=\E\left[I_{kl}^{(1)}\left(Y\right)+I_{kl}^{(2)}\left(Y\right)-I_{kl}^{(3)}\left(Y\right)\right] \,+\,\Delta_n,
\end{equation}
where $\vert  \Delta_n\vert=\mathbf{O}\left(n^{-(1-4c_1-2c_2)}\log n\right)$.
On the other hand,
\begin{multline}\label{i7}
\left|\E\left(I_{kl}^{(2)}\left(Y\right)-I_{kl}^{(3)}\right)\right|
\leq \E\left(\left|\frac{R_{b_n,k}(Y)\widehat{g}_{l,n}(Y)}{f_{b_n}(Y)} - \frac{R_{b_n,k}(Y)R_{b_n,l}(Y)\widehat{f}_{b_n}(Y)}{f_{b_n}(Y)}\right|\right)\\
                                                                     + \E\left(\left|\frac{R_{b_n,l}(Y)\widehat{g}_{k,n}(Y)}{f_{b_n}(Y)} - \frac{R_{b_n,k}(Y)R_{b_n,l}(Y)\widehat{f}_{b_n}(Y)}{f_{b_n}(Y)}\right|\right)
\end{multline}
and
\begin{eqnarray*}
\E\left(\left|\frac{R_{b_n,k}(Y)\widehat{g}_{l,n}(Y)}{f_{b_n}(Y)} - \frac{R_{b_n,k}(Y)R_{b_n,l}(Y)\widehat{f}_{b_n}(Y)}{f_{b_n}(Y)}\right|\right) &=& \E\left(\left|\frac{R_{b_n,k}(Y)}{f_{b_n}(Y)}\right|\left|\widehat{g}_{l,n}(Y) - R_{b_n,l}(Y)\widehat{f}_{b_n}(Y)\right|\right)\nonumber\\
                                                                                                                         &=&
\E\left(\left|\frac{R_{b_n,k}(Y)}{f_{b_n}(Y)}\right|\left|\widehat{g}_{l,n}(Y) - g_l(Y)\frac{\widehat{f}_{b_n}(Y)}{f_{b_n}(Y)}\right|\right)\nonumber\\
                                                                                                                         &\leq&
\E\left(\left|\frac{R_k(Y)}{f(Y)}\right|\left|\widehat{g}_{l,n}(Y) - g_l(Y)\frac{\widehat{f}_{b_n}(Y)}{f_{b_n}(Y)}\right|\right)\label{l84}.
\end{eqnarray*}
Then, using Cauchy-Schwarz inequality and  Lemma \ref{ar22}  we obtain 
\begin{equation*}\label{i21}
\lim_{n\rightarrow +\infty}\E\left(\left|\frac{R_{b_n,k}(Y)\widehat{g}_{l,n}(Y)}{f_{b_n}(Y)} - \frac{R_{b_n,k}(Y)R_{b_n,l}(Y)\widehat{f}_{b_n}(Y)}{f_{b_n}(Y)}\right|\right)= 0.
\end{equation*}
Since, from similar arguments, we also obtain
\begin{equation*}
\lim_{n\rightarrow +\infty}\E\left(\left|\frac{R_{b_n,l}(Y)\widehat{g}_{k,n}(Y)}{f_{b_n}(Y)} - \frac{R_{b_n,k}(Y)R_{b_n,l}(Y)\widehat{f}_{b_n}(Y)}{f_{b_n}(Y)}\right|\right)= 0,
\end{equation*}
we deduce from (\ref{i7}) that  
\begin{equation}\label{lim}
\lim_{n\rightarrow +\infty} \left|\E\left(I_{kl}^{(2)}\left(Y\right)-I_{kl}^{(3)}\right)\right|= 0.
\end{equation}
Moreover, since  $\left|I_{kl,n}^{(1)}(Y)\right|\leq \left|\frac{g_k(Y)g_l(Y)}{f^2(Y)}\right| = \left|R_k(Y)R_l(Y)\right|$, we can apply  the dominated convergence theorem that gives:
\begin{equation}\label{i22}
\lim_{n\rightarrow +\infty}\E\left(I_{kl}^{(1)}(Y)\right)=\E\left(\frac{g_k(Y)g_l(Y)}{f^2(Y)}\right)=\E\Big(R_k(Y)R_l(Y)\Big) =\lambda_{k,l}.
\end{equation}
Then from Eqs. (\ref{i9}), (\ref{lim}) and  (\ref{i22}), it follows
\[
\lim_{n\rightarrow +\infty}\E\left(\widehat{\lambda}^{(n)}_{k,l}\right) =\lambda_{k,l}.
\]
\hfill $\Box$
\subsection{Proof of  Corollary \ref{cor1}}
From Lemma 1 in \cite{ferreyao} (see p. 485)  we have, for any $k\in\{1,\cdots,N\}$,  the inequality  
\[
\Vert\widehat{\eta}_k-\beta_k\Vert_d\leq a_k\Vert\widehat{\Lambda}_n-\Lambda\Vert_\infty,
\] 
where $a_1=2\sqrt{2}/(\lambda_1-\lambda_2)$,  $a_j=2\sqrt{2}/\min(\lambda_{j-1}-\lambda_j,\lambda_{j}-\lambda_{j+1})$ for $j\geq 2$, and $\Vert\cdot\Vert_\infty$ is the matrix norm defined by $\Vert A\Vert_\infty=\sup_{x\in\mathbb{R}^d-\{0\}}\Vert Ax\Vert _d/\Vert x\Vert_d$.  Then, from Theorem \ref{e29} we deduce that  $\widehat{\eta}_k$ converges almost surely to  $\beta_k$ as $n\rightarrow +\infty$. Since $\widehat{\Sigma}_n$ converges almost surely to  $\mathbb{I}_d$  as $n\rightarrow +\infty$, it follows that  $\widehat{\beta}_k$ converges almost surely to  $\beta_k$ as $n\rightarrow +\infty$. 
\hfill $\Box$


\begin{thebibliography}{99}
       \bibitem{Aragon}Aragon~Y, Sarraco~J.   Sliced Inverse Regression (SIR) : an appraisal of small sample alternatives to slicing. Comput.  Statist.  1997; 12:109--130.

 \bibitem{sir2}Duan~N, Li~KC.  Slicing regression: a link-free regression method. Ann. Statist.  1991;  19:505--530.

\bibitem{Bura}Bura~E, Cook~D.   Extending SIR: the weighted chi-square test. J. Amer. Statist. Assoc. 2001; 96:996--1003.

      \bibitem{ferre}Ferr\'e~L.  Determining the dimensionality  in sliced inverse regression and related methods. J. Amer. Statist. Assoc. 1998;  93:132--140.

\bibitem{ferreyao}Ferr\'e~L., Yao~AF. Functional sliced inverse regression analysis. Statistics. 2003; 37:475--488.

     \bibitem{Gine1}Gin\'e~ E,  Guillou~A.  Law of iterated logarithm for censored data. Ann.  Probab. 1999;  27:2042--2067.

     \bibitem{Gine2}Gin\'e~ E,  Guillou~A.  On consistency of kernel density estimators for randomly censored data: rates holding uniformly over adaptive intervals. Ann.    Inst.  Henri Poincar\'e  2001; 37:503--522.  

     \bibitem{Hsing}Hsing~T,  Carroll~ RJ.  An asymptotic theory for sliced inverse regression. Ann. Statist. 1992;  20:1040--1061.

     \bibitem{sir}Li~KC.  Sliced Inverse Regression for Dimension Reduction. J.  Amer.  Statist.  Assoc.1991;  86:316--327.

     \bibitem{Nkiet}Nkiet~GM.  Consistent estimation of the dimensionality in sliced inverse regression. Ann. Inst. Statist. Math.  2008; 60:257--271.

     \bibitem{Rao1}Prakasa Rao~BLS: Nonparametric Functional Estimation.  Orlando: Academic Press; 1983.
\bibitem{saracco}Sarraco~J.  An asymptotic theory for sliced inverse regression. Comm. Stat.- Theory Meth.  1997; 26:2141--2171.

\bibitem{schott}Schott~J.R.   Determining the dimensionality  in sliced inverse regression and related methods. J. Amer. Statist. Assoc. 1998;  89:141--148.

     \bibitem{Talagrand1}Talagrand~M.  Sharper bounds for Gaussian and empirical processes.  Ann. Probab. 1994;  22:28--76.    

     
     \bibitem{tsyb}Tsybakov~AB. Introduction to nonparametric estimation. Paris: Springer; 2009.
   
     \bibitem{vellila}Vellila~S.  Assesing the number of linear components in a general regression problem.   J. Amer. Statist. Assoc. 1998; 93:1088--1098.

     \bibitem{Asymptotics}Zhu~LX,   Fang~KT.  Asymptotics for kernel estimate of sliced inverse regression.   Ann.  Statist  1996; 24:  1053--1068.
    
   \end{thebibliography}
\end{document}